%% file: op4.tex
\documentclass[12pt]{amsart}
\usepackage{amsmath, amssymb, latexsym}

\numberwithin{equation}{section}

\input{defs}

\newtheorem{theorem}{Theorem}[section]
\newtheorem{lemma}[theorem]{Lemma}

\begin{document}
\title[Polynomial Szeg\H o class]
{The Szeg\H o class with a polynomial weight}

\author{S. Denisov}
\author{S. Kupin}

\email{denissov@its.caltech.edu}
\email{kupin@cmi.univ-mrs.fr}


\date{March 15, 2004; preliminary version}

\address{Department of Mathematics 253-37, Caltech, Pasadena CA 91125, USA.}
\address{CMI, Universit\'e de Provence, 39 rue Joliot  Curie, 13453
  Marseille Cedex 13, France.}

\begin{abstract}
Let $p$ be a trigonometric polynomial, nonnegative on the unit
circle $\bt$. We say that a measure $\s$ on $\bt$ belongs to the
polynomial Szeg\H o class, if $d\s(e^{i\th})=\s'_{ac}(e^{i\th})\,
d\th+d\s_s(e^{i\th})$, $\s_s$ is singular, and
$$
\int^{2\pi}_0 p(e^{i\th})\log \s'_{ac}(e^{i\th})\, d\th>-\infty
$$
For the associated orthogonal polynomials $\{\p_n\}$, we obtain
pointwise asymptotics inside the unit disc $\mathbb{D}$. Then we
show that this asymptotics holds in $L^2$-sense on the unit
circle. As a corollary, we get existence of certain modified wave
operators.
\end{abstract}

\maketitle

\vspace{-0.7cm}
\section*{Introduction}

Let $\s$ be a non-trivial Borel probability measure on  the unit
circle \mbox{$\bt=\{z:|z|=1\}$}. Consider orthonormal polynomials
$\{\p_n\}$ with respect to the measure,
$$
\int_\bt\p_n\ovl{\p_m}\, d\s=\d_{nm}
$$
where $\d_{nm}$ is the Kronecker's symbol. Sometimes, it is more convenient to
work with monic orthogonal polynomials $\{\pp_n\},\
\pp_n(z)=z^n+a_{n,n-1}z^{n-1}+\ldots+a_{n,0}$. These polynomials satisfy
$$
\int_\bt\pp_n\ovl{\pp_m}\, d\s=c_n\d_{nm}
$$
with $c_n=||\pp_n||^2_\s=\int_\bt |\pp_n|^2\, d\s$.

It is well-known \cite{ge1,si1} that polynomials $\{\pp_n\}$
generate a sequence \mbox{$\{\a_n\}, |\a_n|<1$},  of the so-called
Verblunsky coefficients through recurrence relations
$$
\begin{array}{lcl}
\pp_{n+1}(z)&=&z\pp_n(z)-\bar\a_n\pp^*_n(z), \pp_0(z)=1\\
\pp^*_{n+1}(z)&=&\pp^*_n(z)-\a_nz\pp_n(z), \pp^*_0(z)=1
\end{array}
$$
where $\pp^*_n(z)=z^n\ovl{\pp_n(1/\bar z)}$. Conversely, the
measure $\s$ (and polynomials $\{\p_n\}$) are completely
determined by the sequence $\{\a_k\}$ of its Verblunsky
coefficients. Hence, it is natural to study the sequence
$\{\a_k\}$ and polynomials $\{\p_n\}$ in terms of $\s$ and vice
versa.

We say that $\s$ is a Szeg\H o measure ($\s\in\ss$, for brevity),
if $d\s=d\s_{ac}+d\s_s=\s'_{ac}dm+d\s_s$ and the density
$\s'_{ac}$ of the absolutely continuous part is such  that
$$
\int_\bt \log \s'_{ac}\, dm>-\infty
$$
Here, the singular part of $\s$ is denoted by $\s_s$, and $m$ is
the probability Lebesgue  measure on $\mathbb{T}$,
$dm(t)=dt/(2\pi it)=1/(2\pi)\, d\th,\ t=e^{i\th}\in\bt$.

The following results are classical.

\begin{theorem}[{\cite{ge1,sz}}]\label{t01}
The following assertions are equivalent
\begin{itemize}
\item[\it i)\ \ ] the sequence $\a$ is in $l^2(\mathbb{Z}_+)$,
\item[\it ii)\ ] the measure $\s$ belongs to the Szeg\H o class,
\item[\it iii)] analytic polynomials are not dense in $L^2(\s)$.
\end{itemize}
\end{theorem}

The last statement of the theorem can be made more precise. Namely, we
have
\begin{equation}
d(\cp_1,0)^2_{L^2(\s)}=\inf_{f\in \cp_1} ||f||^2_\s=\exp \int_\bt
\log\s'_{ac}\, dm \label{e02}
\end{equation}
where $\cp_1$ is the set of analytic polynomials $f$ with the
property $f(0)=1$.

If $\s\in\ss$, we  define a function $D$ in the Hardy class
$H^2(\bd)$ on the unit disk $\bd=\{z:|z|<1\}$ as
\begin{equation}\label{e01}
D(z)=\exp\lt(\frac12 \int_\bt \frac{t+z}{t-z}\, \log\s'_{ac}(t)\,dm(t)\rt)
\end{equation}

\begin{theorem}[{\cite{ge1,sz}}]\label{t02}
Let $\s\in\ss$. Then
$$
\lim_{n\to\infty} D(z)\p^*_n(z)=1
$$
for every $z\in\bd$, and, moreover,
$$
\lim_{n\to\infty} \int_\bt |D\p^*_n-1|^2\, dm=0
$$
that is, the convergence is in the $L^2$-sense on the unit circle.
\end{theorem}

A modern presentation and recent advances in this direction can be found in
\cite{khr1,si1}.

It seems interesting to obtain similar results for
different classes of measures. Consider
a trigonometric polynomial $p$ with the property $p(t)\ge 0,
t\in\bt$. Without loss of generality we can assume it is in the form
\begin{equation}\label{e03}
p(t)=\prod_{k=1}^N |t-\z_k|^{2\kappa_k}
\end{equation}
where $\{\z_k\}$ are points on $\bt$ and $\kappa_k$ are their
``multiplicities''. We say that $\s$ is in the polynomial Szeg\H o
class (i.e., $\s$ is a (pS)-measure or $\s\in\ps$, to be brief),
if $d\s=\s'_{ac}dm+d\s_s$, $\s_s$ being the singular part of the
measure, and
\begin{equation}\label{e04}
\int_\bt p(t)\log\s'_{ac}(t)\, dm(t)>-\infty
\end{equation}
We give counterparts of Theorems \ref{t01} and \ref{t02} for
orthogonal polynomials with respect to (pS)-measures in the next
section. Then we construct modified wave operators for the
corresponding CMV-representations and obtain relations similar to
\eqref{e02}.

\medskip\nt
{\it Acknowledgments.}\ \ The authors are grateful to B.~Simon and
P.~Yuditskii for helpful discussions.

\section{Results}
\label{s1}
\subsection{}
We fix the polynomial $p$ (see \eqref{e03}) for the rest of this
paper. For the sake of transparency, assume $\kappa_k=1$; the
discussion of the general case follows the same lines.

Let $\cc$ and $\cc_0$ be the CMV-representations connected to $\s$
and $m$ (see \mbox{\cite[Ch.~4]{cmv, si1}}), and
$\rk(\cc-\cc_0)<\infty$. Recall that the function $D$ appearing in
\eqref{e01} can be represented as
$$
D(z)=\exp\lt(t_0+\sum^\infty_{k=1}\frac{\tr(\bar\cc^k-\bar\cc^k_0)}{k}\,
z^k\rt)
$$
Here,
$$
t_0=\sum_k\log\r_k=\sum_k\log(1-|\a_k|^2)^{1/2}
$$
and $\{\a_k\}$ are the Verblunsky coefficients corresponding to $\s$.

Furthermore, the polynomial $p$ \eqref{e03}
defines an analytic polynomial $P$ via the formulas
$$
p_1=2P_+p,\quad P'(t)=\frac{p_1(t)-p_1(0)}{t}
$$
where $P_+:L^2(\bt)\to H^2(\bd)$ is the Riesz projector \cite[Ch. 3]{ga}.

\begin{lemma}\label{l1}
Let $\rk(\cc-\cc_0)<\infty$. Then
\begin{equation}\label{e1}
\int_\bt p\log\s'_{ac}\, dm=a_0t_0+\re\tr (P(\cc)-P(\cc_0))
\end{equation}
where
$$
a_0=p_1(0)=2\int_\bt p\, dm
$$
\end{lemma}

We denote the right-hand side of equality \eqref{e1} by
$\pps(\cc)$ and we rewrite it in a slightly different form. To
this end, we consider the shift $S:l^2(\bz_+)\to l^2(\bz_+)$,
given by $Se_k=e_{k+1}$. For a bounded  operator $A$ on
$l^2(\bz_+)$, we look at $\t(A)=S^*AS$. In particular, the matrix
of $\t^k(A), k\in\bz_+$, is obtained from the matrix of $A$ by
dropping the first k rows and columns. Going back to \eqref{e1},
we get that
\begin{eqnarray*}
\pps(\cc)&=&\sum^{\infty}_{k=0}\lt\{a_0\log\r_k +
\re ((P(\cc)-P(\cc_0))e_k,e_k)\rt\}\\
&=&\sum^{2N+1}_{k=0}\lt\{a_0\log\r_k+\re ((P(\cc)-P(\cc_0))e_k,e_k)\rt\}
+\sum^\infty_{k=2N+2}\psi\circ\t^k(\cc)
\end{eqnarray*}
where $\psi(\cc)=a_0\log\r_{2N+2}+\re((P(\cc)-P(\cc_0))e_{2N+2}, e_{2N+2})$.

Repeating word-by-word arguments from \cite{nvyu}, Lemma 3.1, we see
that there exists a function $\g$, depending on $l=4N+5$ arguments,
such that
$$
\psi(x_1,\ldots,
x_l)=\eta(x_1,\ldots,x_l)-\g(x_2,\ldots,x_l)+\g(x_1,\ldots,x_{l-1})
$$
and $\eta(x_1,\ldots, x_l)\le 0$ for any collection
$(x_1,\ldots,x_l)$.
Now, we put
\begin{eqnarray*}
\pp(\cc)&=&\int_\bt p(t)\log\s'_{ac}(t)\, dm(t) \\
\ti\pps(\cc)&=&\sum^{2N+1}_{k=0}\big\{a_0\log\r_k
+\re ((P(\cc)-P(\cc_0))e_k,e_k)\big\} \\
&+&\sum^\infty_{k=2N+2}\eta\circ\t^k(\cc)+\g\circ\t^{2N+2}(\cc)
\end{eqnarray*}
Consequently, for a CMV-representation $\cc$ having
$\rk(\cc-\cc_0)<\infty$, equality \eqref{e1} reads as
$$
\pp(\cc)=\pps(\cc)=\ti\pps(\cc)
$$
The following theorem holds.

\begin{theorem}[{\cite[Theorem~1.4]{nvyu}}]\label{t1}
A measure $\s$ is polynomially Szeg\H o (see \eqref{e04}) if and
only if \ $\ti\pps(\cc)>-\infty$. Moreover, in this case \
$\pp(\cc)=\ti\pps(\cc)$.
\end{theorem}

\subsection{}
We turn now to the description of asymptotical properties of
orthogonal polynomials for (pS)-measures. Consider a modified
Schwarz kernel
$$
K(t,z)=\frac{t+z}{t-z}\,\frac{q(t)}{q(z)}
$$
where $q(t)=C(\prod_k(t-\z_k)^2)/t^N$, and the constant $C,|C|=1,$ is
chosen in a way that $q(t)\in\br$ for $t\in\bt$ (i.e.,
$C=(\prod_k(-\z_k))^{-1}$). Furthermore, define
\begin{eqnarray*}
\ti D(z)&=&\exp\lt(\frac12\int_\bt K(t,z)\log\s'_{ac}(t)\, dm(t)\rt)\\
\ti \p^*_n(z)&=&\exp\lt(\int_\bt K(t,z)\log|\p^*_n(t)|\, dm(t)\rt)
\end{eqnarray*}
The functions $\{\ti\p^*_n\}$ are called (reversed) modified
orthogonal polynomials with respect to $\s$. It can be readily
seen that $|\ti D|^2=\s'_{ac}$ and $|\ti\p^*_n|=|\p^*_n|=|\p_n|$
a.e. on $\bt$. Furthermore, we see that $\ti\p^*_n=\psi_n\p^*_n$,
where
\begin{eqnarray}
\psi_n(z)&=&\exp\lt(\int_\bt\frac{t+z}{t-z}\lt(\frac{q(t)}{q(z)}-
1\rt)\log|\p^*_n(t)|\, dm(t)\rt)\nonumber\\
&=&\exp\lt(A_{0n}+\sum^N_{k=1}\lt(A_{kn}\frac{z+\z_k}{z-\z_k}
+B_{kn}\lt\{\frac{z+\z_k}{z-\z_k}\rt\}^2\rt)\rt)\label{e2}
\end{eqnarray}
and $A_{0n},B_{kn}\in i\br$, $A_{kn}\in\br$. The coefficients
$\{A_{0n},A_{kn},B_{kn}\}_{k,n}$ can be expressed in a closed form
through $\{\a_k\}$.

\begin{theorem}\label{t2}
Let $\s\in\ps$. Then
\begin{itemize}
\item[\it i) \ ] for any $z\in\bd$,
$$
\lim_{n\to\infty} \ti D(z)\ti\p^*_n(z)=\lim_{n\to\infty}\ti
D(z)\psi_n(z)\p^*_n(z)=1
$$
\item[\it ii)] we also have
$$
\lim_{n\to\infty} \int_\bt|\ti D\ti\p^*_n-1|^2\, dm=0
$$
\end{itemize}
\end{theorem}

The proof of the first claim of the theorem is partially based on
the sum
rules proved in Theorem \ref{t1}. One of the main facts leading to the
second claim is that
$$
\lim_{n\to\infty}\int_I|\ti D\ti\p^*_n-1|^2\, dm=0
$$
for any closed arc $I\subset \bt$ that does not contain points
$\{\z_k\}$. We prove the latter relation observing that, for small
$\ep>0$,
$$
|\ti D\ti\p^*_n(z)|\le\frac{C(\ep)}{\sqrt{1-|z|}}
$$
for $z\in \bd\bsl(\cup_k B_\ep(\z_k))$, $B_\ep(\z)=\{z:
|z-\z|<\ep\}$.

\subsection{}
We now use asymptotics described in the last subsection, to construct
modified wave operators.

Let $\cf_0:L^2(m)\to l^2(\bz_+), \cf:L^2(\s)\to l^2(\bz_+)$ be the
Fourier transforms associated to the CMV-representations $\cc$ and
$\cc_0$, see \cite[Ch.~4]{si1}. Recall that
$$
\cc=\cf z\cf^{-1},\quad \cc_0=\cf_0 z\cf^{-1}_0
$$

\begin{theorem}\label{t3}
Let $\s\in\ps$. The limits
\begin{equation}\label{e4}
\ti\oo_\pm=\sl_{n\to\pm\infty}\, e^{W(2n,\cc)}\cc^n\cc_0^{-n}
\end{equation}
exist. Here
$$
W(C,n)=A_{0n}+\sum^N_{k=1}\lt(A_{kn}\frac{\cc+\z_k}{\cc-\z_k}
+B_{kn}\lt\{\frac{\cc+\z_k}{\cc-\z_k}\rt\}^2\rt)
$$
and coefficients $\{A_{0n},A_{kn},B_{kn}\}$ are defined in
\eqref{e2}. We also have
$$
\cf^{-1}\ti\oo_+\cf_0=\chi_{E_{ac}}\,\frac1{\ti D},\quad
\cf^{-1}\ti\oo_-\cf_0=\chi_{E_{ac}}\,\frac1{\bar{\ti D}}
$$
where $E_{ac}=\bt\bsl\mathrm{supp}\, \s_s$.
\end{theorem}
As a simple corollary, we get the existence of modified wave
operators for the pair ($\cc_0, \cc$). The natural (and open)
question is to prove modified wave operators for the pair ($\cc,
\cc_0$).

\subsection{}
We now briefly discuss a variational problem related to
(pS)-measures. We put $\cp'_0$ to be the set of polynomials $g$
analytic on $\bd$ with the property $g\not =0$ on $\bd$ and
$g(0)>0$. Furthermore, for a $g\in \cp'_0$, we set
$$
\l(g)=\exp\lt(\int_\bt p\log|g|\, dm\rt)
$$
and define
$\cp'_1=\{g: g\in\cp'_0,\ \l(g)=1\}$.

\begin{theorem}\label{t4}
Let $d\s=w\, dm+d\s_s$. Then
\begin{eqnarray}
\exp\lt(\int_\bt p\log\frac wp\,dm\rt)&\le&\inf_{g\in\cp'_1}||g||^2_\s=
\inf_{
\dsp
\begin{array}{c}
g\in\cp'_0,\\
||g||_\s\le 1
\end{array}}
\frac1{|\l(g)|^2} \label{e6}\\
&\le&\exp\lt(\int_\bt p\log w\, dm\rt) \nonumber
\end{eqnarray}
\end{theorem}

Remind that $\s$ is a Szeg\H o measure if and only if the
system $\{e^{iks}\}_{k\in\bz}$ is uniformly minimal in
$L^2(\s)$. Saying that $\s$ is a (pS)-measure translates into the
uniform minimality of another system, $\{e^{ik\nu(s)}\}_{k\in\bz}$,
in the same
space $L^2(\s)$. Above,
$$
\nu(s)=C_0\int^s_0 p(e^{is'})\, ds'
$$
where $s,s'\in [0,2\pi]$ and the constant $C_0$ comes from the
condition $C_0\int_\bt p\, dm=1$, see \cite{nvyu}, Lemma 2.2. It
seems to be an interesting question to translate these results to
the language of Gaussian processes \cite{dm}.

It was proved recently in \cite[Ch. 2]{si1} that
$\s\in\mathrm{(p}_0\mathrm{S)}$ with
$$
p_0(t)=\frac12\, |1-t|^2=1-\cos\th
$$
if and only if $\{\a_k\}\in l^4$ and $\{\a_{k+1}-\a_k\}\in l^2$
(above, $t=e^{i\th}$).
Theorems \ref{t1}--\ref{t4} readily apply to this special case.
In particular, we have
\begin{eqnarray*}
K_0(t,z)&=&\frac{t+z}{t-z}\, \frac{(t-1)^2}{t}\, \frac{z}{(1-z)^2}\\
\ti D_0(z)&=&\exp\lt(\frac12\int_\bt K_0(t,z)\log\s'_{ac}(t)\, dm(t)\rt)
\end{eqnarray*}
and
$$
\psi_n(z)=\exp\lt(A_n\lt\{\frac{1+z}{1-z}-1\rt\} +
B_n\lt\{\lt(\frac{1+z}{1-z}\rt)^2-1\rt\}\rt)
$$
where
$$
A_n=\sum^n_{k=0}\log (1-|\a_k|^2)^{1/2},\quad
B_n=\frac i4\, \im\Big(\a_0-\sum^n_{k=1}\bar\a_{k-1}\a_k\Big)
$$

\end{document}

%% file: defs.tex
\newcommand{\bd}{\mathbb{D}}

\newcommand{\br}{\mathbb{R}}

\newcommand{\bt}{\mathbb{T}}

\newcommand{\bz}{\mathbb{Z}}

\renewcommand{\a}{\alpha}

\renewcommand{\l}{\lambda}

\newcommand{\s}{\sigma}

\renewcommand{\r}{\rho}

\newcommand{\p}{\varphi}

\newcommand{\pp}{\Phi}

\newcommand{\pps}{\Psi}
\renewcommand{\t}{\tau}

\renewcommand{\th}{\theta}

\renewcommand{\d}{\delta}

\newcommand{\oo}{\Omega}

\newcommand{\g}{\gamma}

\newcommand{\ep}{\varepsilon}

\newcommand{\z}{\zeta}

\newcommand{\nt}{\noindent}

\newcommand{\bsl}{\backslash}

\newcommand{\ovl}{\overline}

\newcommand{\ti}{\tilde}

\newcommand{\lt}{\left}

\newcommand{\rt}{\right}

\newcommand{\dsp}{\displaystyle}

\renewcommand{\ss}{\mathrm{(S)}}

\newcommand{\ps}{\mathrm{(pS)}}
\newcommand{\cf}{\mathcal{F}}

\newcommand{\cc}{\mathcal{C}}

\newcommand{\cp}{\mathcal{P}}

\newcommand{\tr}{\mathrm{tr}\,}

\newcommand{\re}{\mathrm{Re}\,}

\newcommand{\im}{\mathrm{Im}\,}

\newcommand{\rk}{\mathrm{rank}\,}

\renewcommand{\sl}{\mathrm{s-lim}\,}